\documentclass[12pt,letterpaper]{article}
\usepackage[english]{babel}
\usepackage{hyperref}

\date{}
\newcommand\proofr[1][\hspace{-1em}]{\par{\em Proof\hspace{1em}#1:~}}

\def\proofend{$\bigtriangleup$\vspace{0.3em}\par}

\newcounter{theorem}
\newenvironment{theorem}[1][\hspace{-1.0ex}]%
 {\par\addvspace{2mm}\indent\refstepcounter{theorem} \textbf{Theorem\hspace{1.0ex}{\rm#1}.~}\sl}%
 {\par\addvspace{2mm}\rm}
\newcounter{remark}
\newenvironment{remark}[1][\hspace{-1.0ex}]%
 {\par\addvspace{2mm}\indent\refstepcounter{remark} \textbf{Remark~\theremark\hspace{1.0ex}{\rm#1}.~}}%
 {\par\addvspace{2mm}\rm}
\newcounter{lemma}0
\newenvironment{lemma}[1][\hspace{-1.0ex}]%
 {\par\addvspace{2mm}\indent\refstepcounter{lemma} \textbf{Lemma~\thelemma\hspace{1.0ex}{\rm#1}.~}\sl}%
 {\par\addvspace{2mm}\rm}

\providecommand\href[2]{#2}
\author{Sergey V. Avgustinovich, Denis S. Krotov}

\def\zeroext#1{\bar #1}
\begin{document}


\begin{center}
{\LARGE Embedding in a perfect code%
\footnote{
This is the peer reviewed version of the following article: 
``Embedding in a perfect code'', Journal of Combinatorial Designs 17(5) 2009, 419--423, 
which has been published in final form at 
\url{http://dx.doi.org/10.1002/jcd.20207}.
This article may be used for non-commercial purposes in accordance with 
Wiley Terms and Conditions for Self-Archiving.
\\
The results of the paper were presented at the Workshop
``Coding Theory Days in St. Petersburg'', Oct. 2008, St. Petersburg, Russia.
\\
The authors are with the Sobolev Institute of Mathematics, Novosibirsk, Russia.
E-mail: \{avgust,krotov\}@math.nsc.ru
\\
The research was partially supported by the RFBR grant 07-01-00248 and 08-01-00673}
}
\end{center}

\begin{center}
{\large
{Sergey~V.~Avgustinovich},
{Denis~S.~Krotov}}
\end{center}

\renewcommand\abstractname{}
\begin{abstract}
\begin{center}
A binary $1$-error-correcting code can always be embedded \\in a $1$-perfect code
of some larger length.\end{center}
\end{abstract}

For any  $1$-error-correcting binary code $C$ of length $m$ we will construct a $1$-perfect
binary code $P(C)$ of length $n=2^m-1$ such that fixing the last $n-m$ coordinates by zeroes in $P(C)$
gives $C$.

In particular, any complete or partial Steiner triple system (or any other system
that forms a $1$-code)
can always be embedded in a $1$-perfect code of some length
(compare with \cite{OstPot2007}).
Since the weight-$3$ words of a $1$-perfect code $P$ with $0^n \in P$ form a Steiner triple system,
and the weight-$4$ words of an extended $1$-perfect code $P$ with $0^n \in P$
form a Steiner quadruple system, we have, as corollaries, the following well-known facts:
a patrial Steiner triple (quadruple) system can always be embedded
in a Steiner triple (quadruple) system \cite{Treash:TStoSTS} (\cite{Ganter:QStoSQS})
(these results, as well as many other embedding theorems for Steiner systems,
can be found in \cite{Linder:survey-embed,ColRos}).

Notation:
\begin{itemize}
\item $F^m$ denotes the set of binary $m$-tuples, or binary $m$-words.
\item $\dot F^m := F^m \setminus \{0^m \}$, where $0^m$ is the all-zeroes $m$-word.
\item$F^m$ is considered as a vector space over $GF(2)$ with calculations modulo $2$.
\item${{\Pi}}=\{\pi^{(1)},\ldots,\pi^{(m)}\}=\{(10..0),\ldots,(0..01)\}$ is the natural basis in $F^m$.
\item$n:=2^m-1$.
\item The elements of $F^m$ will be denoted by Greek letters.
\item The elements of $F^n$ will be denoted by overlined letters,
their coordinates being indexed by the elements of $\dot F^m$, e.g.,
$\bar w = \{w_\iota\}_{\iota\in \dot F^m}$; we assume
that the first $m$ coordinates have the indexes $\pi^{(1)},\ldots,\pi^{(m)}$, and
the order of the other $n-m$ indexes does not matter (but fixed).
\item
$\{\bar e^{(\iota)}\}_{\iota\in \dot F^m}$ is the natural basis in $F^n$;
note that $\bar e^{(\pi^{(\iota)})}=(\pi^{(\iota)},0^{n-m})$.
\item
For any $\alpha=(\alpha_1,...,\alpha_m)\in F^m$ denote
$$ \zeroext{\alpha} := (\alpha,0^{n-m});$$
it also holds $ \zeroext{\alpha} =\sum_{i=1}^{m} \alpha_i \bar e^{(\pi^{(i)})}. $
\item
$d(\cdot,\cdot)$ denotes the Hamming \emph{distance} between two words in $F^m$ or $F^n$
(the number of positions in which the words differ).
\item
$<\ldots>$ denotes the linear span of the vectors or sets of vectors between the angle brackets.
\item
The \emph{neighborhood} $\Omega(M)$ of a set $M\subset F^n$ is the set of vectors
at distance at most $1$ from $M$.
\item
A set $C\subset F^m$ is called a \emph{$1$-code} if the neighborhoods of the codewords are disjoint.
\item
A $1$-code $P\subset F^n$ is called a \emph{$1$-perfect code} if $\Omega(P)=F^n$;
in this case, $|P|=2^n/(n+1)$.
\item
The
\emph{Hamming code} $H$ defined as
\begin{equation}\label{eq:H}
H := \{ \bar c\in\{0,1\}^n | \sum_{\alpha\in\dot F^m} c_\alpha \alpha = 0^m   \}
\end{equation}
is a linear $1$-perfect code.
\item
For any $\iota$ from $\dot F^m$ the \emph{linear $\iota$-component} of $H$ is defined as
$$
R_\iota := \{ \bar c\in H |  c_\alpha = c_{\alpha+\iota} \mbox{ for all }  \alpha \in F^m\setminus <\iota>\}
$$
\end{itemize}
(note that $R_\iota$ is a linear subcode of $H$, for all $\iota$).
Since \cite{Vas:nongroup_perfect}, linear components are used for constructing
non-linear $1$-perfect codes.
For the first time, the method of synchronous switching nonintersecting
linear $i$-components with different $i$,
which is exploited in this paper
(to follow our notations, we replace $i$ by Greek letters),
was used in \cite{EV:94}.
This method was also used by different authors for constructing $1$-perfect codes with
specified properties,
such as different ranks and(or) kernels
\cite{EV:94,PheLeV:1995,PheVil:RankKernel,ASH:RankKernel},
trivial automorphism group \cite{AvgSol:1998IT:trivAut,Mal:1998OC:trivAut},
nonsystematicity \cite{AvgSol:1996:nonsyst, PheLeV:1999, Mal:2001:nonsyst}
(to read more about these and other results in the theory of $1$-perfect codes,
see \cite{Hed:2008:survey,Sol:2008:survey}).
The term ``switching'' was introduced in
\cite{PheLeV:1999}; probably, it came from the theory of Steiner triple systems.
In general, this term refers to replacing some part of a $1$-perfect code by some set with
the same cardinality and same neighborhood; the resulting code will also be $1$-perfect.
With this approach, a large number of $1$-perfect codes can be produced \cite{KroAvg:2008:lb}.
$1$-Codes that admit replacing by other $1$-codes with the same neighborhood
deserve the independent study; they exist in $F^n$ for every odd $n$ \cite{VAK:2008}.

%
Being such a $1$-code is the main property of $R_\iota$.
Since our definition of linear components differs from others,
we should prove this fact
(in essence, the following lemma coincides with \cite[Corollary 3.4]{EV:94}).

\begin{lemma}\label{l:O} For any $\bar z$ from $F^n$ it holds that
  $\Omega(R_\iota+\bar z) = \Omega(R_\iota+\bar z+\bar e^{(\iota)}).$
\end{lemma}
\proofr Without loss of generality, assume $\bar z = 0^n$.
Denote $\bar e^{(0^m)}:=0^n$. Then, we have
$$
\Omega(R_\iota)
=\bigcup_{\kappa\in F^m}(R_\iota + \bar e^{(\kappa)})
=\bigcup_{\kappa\in F^m}(R_\iota + \bar e^{(\iota)} + \bar e^{(\kappa+\iota)})
=\bigcup_{\lambda\in F^m}( (R_\iota + \bar e^{(\iota)}) + \bar e^{(\lambda)})
=\Omega(R_\iota + \bar e^{(\iota)})
$$
because $\bar e^{(\iota)} + \bar e^{(\kappa)} + \bar e^{(\kappa+\iota)} \in R_\iota$ for all $\kappa\in F^m$.
\proofend

\begin{lemma}\label{l:RR} Every element $\bar c$ of $<R_\iota,R_\kappa>$ satisfies
\begin{equation}\label{eq:RR}
 c_{\alpha} + c_{\alpha+\iota}+c_{\alpha+\kappa} + c_{\alpha+\iota+\kappa}=0
\mbox{ for all } \alpha \in F^m\setminus <\iota,\kappa>
\end{equation}
\end{lemma}
\proofr
By the definition, the elements of $R_\iota$ and $R_\kappa$ satisfy (\ref{eq:RR}).
Thus, the elements of their linear span also satisfy (\ref{eq:RR}).
\proofend

The following lemma is the crucial part of our reasoning.

\begin{lemma}\label{l:disj}
  For any $\iota,\kappa\in\dot F^m$ at distance at least $3$ from $0^m$
  and from each other,
  the $\iota$-component
  $R_\iota+ \zeroext{\iota} + \bar e^{(\iota)}$
  and
  the $\kappa$-component
  $R_\kappa+ \zeroext{\kappa} + \bar e^{(\kappa)}$
  are disjoint
  and do not contain $0^n$.
\end{lemma}
\proofr
By general algebraic reasons, it is enough to show that
$\bar w:= \zeroext{\iota} + \bar e^{(\iota)} + \zeroext{\kappa} + \bar e^{(\kappa)}$,
 does not belong to
$<R_\iota,R_\kappa>$.
Let $j$ be a nonzero coordinate of $\iota+\kappa$.
Then $\pi^{(j)}$ is the index of
a nonzero coordinate of $\bar w$;
the indexes of the other nonzero coordinates
also belong to ${\Pi} \cup \{\iota,\kappa\}$.
But, since the mutual distances between $0^m$, $\iota$, $\kappa$, and $\iota+\kappa$
are not less than $3$,
the indexes $\pi^{(j)}+\iota$, $\pi^{(j)}+\kappa$, $\pi^{(j)}+\iota+\kappa$
do not belong to ${\Pi} \cup \{\iota,\kappa\}$.
So, we have
$w_{\pi^{(j)}}+w_{\pi^{(j)}+\iota}+w_{\pi^{(j)}+\kappa}+w_{\pi^{(j)}+\iota+\kappa} = 1+0+0+0=1$,
and, by Lemma~\ref{l:RR}, $\bar w \not\in <R_\iota,R_\kappa>$.

By a similar argument, neither $R_\iota+ \zeroext{\iota} + \bar e^{(\iota)}$
nor $R_\kappa+ \zeroext{\kappa} + \bar e^{(\kappa)}$ contains $0^m$.
\proofend

\begin{theorem}
  Let $C\subset F^m$ be a $1$-code containing $0^m$; put $\dot C := C \setminus \{0^m\}$.
  Then the set
 $$ P(C) := \left( H \setminus \bigcup_{\iota\in \dot C}(R_\iota+ \zeroext{\iota} + \bar e^{(\iota)})
 \right) \cup \bigcup_{\iota\in \dot C}(R_\iota+ \zeroext{\iota})$$
 is a $1$-perfect code in $F^n$. Moreover,
\begin{equation}\label{eq:CtoP}
  C= \{ \iota \in F^m | (\iota,0^{n-m})\in P(C)\}.
\end{equation}
\end{theorem}
\proofr
We note that, by (\ref{eq:H}),  $\zeroext{\iota} + \bar e^{(\iota)}$
belongs to $H$ for all $\iota$; so, $R_\iota+ \zeroext{\iota} + \bar e^{(\iota)} \subset H$ for all $\iota$.

By Lemma~\ref{l:disj}, the sets $R_\iota+ \zeroext{\iota} + \bar e^{(\iota)}$, $\iota\in \dot C$,
are mutually disjoint.
Since they are subsets of a $1$-perfect code, their neighborhoods are also mutually
disjoint. Taking into account Lemma~\ref{l:O}, we see that $P(C)$ is a $1$-perfect code
by the definition.

It is easy to see that\\
(*) \emph{the only word in $H$ that has the form $(\alpha,0^{n-m})$
is the all-zeroes word.}

Furthermore, only $\zeroext{\iota}$ has such form in $R_\iota+ \zeroext{\iota}$, i.e.,\\
(**) \emph{if for some $\kappa\in F^m$
we have $(\kappa,0^{n-m}) \in R_\iota+ \zeroext{\iota}$,
then $\kappa = \iota$.}
Indeed, assume that $\zeroext{\kappa}=(\kappa,0^{n-m})\in R_\iota+ \zeroext{\iota}$.
Then $\zeroext{\kappa} + \zeroext{\iota} \in R_\iota \subset H$.
By (*), we have $\zeroext{\kappa} + \zeroext{\iota}=0^n$,
which proves the claim (**).

From (*) and (**) we conclude that (\ref{eq:CtoP}) is implied by the definition of $P(C)$.
\proofend

\begin{remark}
It is not difficult to see that ${\mathrm{rank}}(P(C))={\mathrm{rank}}(H)+{\mathrm{rank}}(C)$, where the rank of a code is the dimension
  of its linear span.
  So, if $m$ is sufficiently large, we can construct a code $P(C)$ of any rank from $2^m-m$ to $2^m-1$.
 This gives  (for large lengths) an interesting treatment of the result \cite[Proposition 6.2]{EV:94}
 on the existence of $1$-perfect
 codes of any ${\mathrm{rank}}\geq {\mathrm{rank}}(H)$.
\end{remark}
\begin{remark} As proved in \cite{Mal:2001:nonsyst}, it is enough to switch only
  seven disjoint cosets $R_\iota+z_\iota$,  $R_\kappa+z_\kappa$,  $R_\lambda+z_\lambda$,
  $R_{\iota+\kappa}+z_{\iota+\kappa}$, $R_{\iota+\lambda}+z_{\iota+\lambda}$,
  $R_{\kappa+\lambda}+z_{\kappa+\lambda}$, $R_{\iota+\kappa+\lambda}+z_{\iota+\kappa+\lambda}$
  in $H$ for some linearly independent $\iota$, $\kappa$, $\lambda$ and some
  $z_\iota,z_\kappa,z_\lambda,z_{\iota+\kappa},z_{\iota+\lambda},z_{\kappa+\lambda},z_{\iota+\kappa+\lambda}\in H$
  to obtain a nonsystematic $1$-perfect code. So, starting from $m=6$, we
  can construct a nonsystematic $1$-perfect code $P(C)$ with $|C|=8$.
\end{remark}

It would be interesting to define the minimal value of $n(m)$ such that any $1$-code of length $m$
can always be embedded in a $1$-perfect code of length $n(m)$. Recently, the similar question for
Steiner triple systems was solved \cite{BryHor:Lindner}:
any partial Steiner triple system of order $u$ can always be embedded in a Steiner triple systems
of order $v$ if and only if  $v\geq 2u+1$ and $v \equiv 1,3 \bmod 6$
(the last is the necessary and sufficient condition for the existence
of Steiner triple systems of order $v$). Taking into account the relation between
perfect codes and Steiner triple systems discussed in the introduction, we can state
that $n(m)\geq 2m+1$; the Theorem gives $n(m)\leq 2^m-1$.

The authors would like to thank the referees for their useful comments.

\providecommand\href[2]{#2} \providecommand\url[1]{\href{#1}{#1}}
  \def\DOI#1{{\small {DOI}:
  \href{http://dx.doi.org/#1}{#1}}}\def\DOIURL#1#2{{\small{DOI}:
  \href{http://dx.doi.org/#2}{#1}}}\providecommand\bbljun{June}

\end{document}



\begin{center}
{\LARGE Вложение в совершенный код%
\footnote{Настоящий текст является авторским переводом статьи ``Embedding in a perfect code'', Journal of Combinatorial Designs 17(5) 2009, 419--423, DOI \href{http://dx.doi.org/10.1002/jcd.20207}{10.1002/jcd.20207}. \copyright~2009 Wiley Periodicals, Inc.
\\ \indent
Результаты работы были доложены на международной конференции
``Coding Theory Days in St. Petersburg'', 6--10 октября 2008\,г., Санкт"=Петербург, Россия.
\\ \indent
Адрес авторов: Институт математики им. С. Л. Соболева СО РАН, проспект Академика Коптюга 4, Новосибирск 630090, Россия.
E-mail: \{avgust,krotov\}@math.nsc.ru
\\ \indent
Работа выполнена при финансовой поддержке РФФИ (проекты 07-01-00248 и 08-01-00673).
}}
\end{center}

\begin{center}
{\large
{С.~В.~Августинович},
{Д.~С.~Кротов}}
\end{center}

\renewcommand\abstractname{}
\begin{abstract}
Любой двоичный код, исправляющий $1$ ошибку, всегда можно вложить в $1$-совершенный код некоторой большей длины.
\end{abstract}

Для любого двоичного кода $C$ длины $m$, исправляющего $1$ ошибку, мы построим $1$-совершенный двоичный код $P(C)$ длины $n=2^m-1$ такой, что при фиксации последних $n-m$ координат нулями из него получается код $C$.

В частности, любая полная или частичная система троек Штейнера (или любая другая система, которая образует код, исправляющий $1$ ошибку)
всегда может быть вложена в $1$-совершенный код некоторой длины (см. для сравнения \cite{ru:OstPot2007}).
Поскольку слова веса $3$ в $1$-совершенном коде $P$ таком, что $0^n \in P$, образуют систему троек Штейнера, и слова веса $4$ расширенного $1$-совершенного кода $\overline P$ такого, что $0^n \in \overline P$, образуют систему четверок Штейнера, мы имеем в качестве следствия следующие хорошо известные факты: частичная система троек (четверок) Штейнера всегда вложима в систему троек (четверок) Штейнера \cite{ru:Treash:TStoSTS} (\cite{ru:Ganter:QStoSQS})
(эти результаты, как и много других теорем о вложениях систем Штейнера, можно найти в \cite{ru:Linder:survey-embed,ru:ColRos}).

Введем необходимые обозначения:
\begin{itemize}
\item $F^m$ "--- множество двоичных наборов длины $m$, или двоичных $m$-слов.
\item $\dot F^m := F^m \setminus \{0^m \}$, где $0^m$ есть $m$-слово из всех нулей.
\item$F^m$ рассматривается как векторное пространство над полем $GF(2)$ с вычислениями по модулю $2$.
\item${{\Pi}}=\{\pi^{(1)},\ldots,\pi^{(m)}\}=\{(10..0),\ldots,(0..01)\}$ "--- естественный базис в $F^m$.
\item$n:=2^m-1$.
\item Элементы множества $F^m$ будем обозначать греческими буквами.
\item Элементы множества $F^n$ будем обозначать буквами с чертой сверху, их координаты будем индексировать элементами множества $\dot F^m$, т.е.
$\bar w = \{w_\iota\}_{\iota\in \dot F^m}$; договоримся, что первые $m$ координат имеют индексы $\pi^{(1)},\ldots,\pi^{(m)}$, порядок остальных $n-m$ индексов произвольный (но считается фиксированным).
\item
$\{\bar e^{(\iota)}\}_{\iota\in \dot F^m}$ "--- естественный базис в $F^n$;
заметим, что $\bar e^{(\pi^{(\iota)})}=(\pi^{(\iota)},0^{n-m})$.
\item
Для любого $\alpha=(\alpha_1,...,\alpha_m)\in F^m$ обозначим
$$ \zeroext{\alpha} := (\alpha,0^{n-m});$$
выполняется соотношение $ \zeroext{\alpha} =\sum_{i=1}^{m} \alpha_i \bar e^{(\pi^{(i)})}. $
\item
$d(\cdot,\cdot)$ обозначает расстояние \emph{Хэмминга} между двумя словами в $F^m$ или $F^n$ (число позиций, в которых слова различаются).
\item
$<\ldots>$ обозначает линейную оболочку векторов или множеств векторов, записанных внутри угловых скобок.
\item
 \emph{Окрестность} $\Omega(M)$ множества $M\subset F^n$ есть множество векторов на расстоянии не более $1$ от $M$.
\item
Множество $C\subset F^m$ называется \emph{$1$-кодом}, если окрестности кодовых слов попарно не пересекаются.
\item
$1$-Код $P\subset F^n$ называется \emph{$1$-совершенным кодом}, если $\Omega(P)=F^n$;
в этом случае $|P|=2^n/(n+1)$.
\item
\emph{Код Хэмминга} $H$, определяемый соотношением
\begin{equation}\label{eq:H}
H := \{ \bar c\in\{0,1\}^n | \sum_{\alpha\in\dot F^m} c_\alpha \alpha = 0^m   \},
\end{equation}
есть линейный $1$-совершенный код.
\item
Для любого $\iota$ из $\dot F^m$ \emph{линейная $\iota$-компонента} кода $H$ определяется следующим образом:
$$
R_\iota := \{ \bar c\in H |  c_\alpha = c_{\alpha+\iota} \mbox{ для всех }  \alpha \in F^m\setminus <\iota>\}
$$
\end{itemize}
(заметим, что для любого $\iota$ множество $R_\iota$ является линейным подкодом кода $H$).
Начиная с работы \cite{ru:Vas:nongroup_perfect.ru}, линейные компоненты используются для построения нелинейных $1$-совершенных кодов.
Впервые метод одновременного свитчинга (сдвига) непересекающихся линейных $i$-компонент с различными $i$, который используется и в данной работе (чтобы следовать нашим обозначениям, мы заменяем $i$ греческими буквами), был применен в работе \cite{ru:EV:94}.
Этот метод использовался различными авторами для построения $1$-совершенных кодов с особыми свойствами, такими как различные ранги и(или) ядра \cite{ru:EV:94,ru:PheLeV:1995,ru:PheVil:RankKernel,ru:ASH:RankKernel},
тривиальная группа автоморфизмов \cite{ru:AvgSol:1998IT:trivAut,ru:Mal:1998OC:trivAut},
несистематичность \cite{ru:AvgSol:1996:nonsyst.ru, ru:PheLeV:1999, ru:Mal:2001:nonsyst.ru}
(для ознакомления с этими и некоторыми другими результатами теории $1$-совершенных кодов см. обзоры \cite{ru:Hed:2008:survey,ru:Sol:2008:survey}).
Термин <<свитчинг>> был предложен в работе \cite{ru:PheLeV:1999}; вероятно он берет начало из теории систем троек Штейнера.
В самой общей интерпретации, этот термин означает замену некоторой части $1$-совершенного кода некоторым множеством той же мощности и с той же окрестностью, при этом получившийся код также будет $1$-совершенным.
При помощи такого подхода можно построить большое число $1$-совершенных кодов \cite{ru:KroAvg:2008:lb}.
$1$-Коды, которые допускают замену на другие $1$-коды с той же окрестностью, заслуживают независимого изучения, они существуют в $F^n$ при каждом нечетном $n$ \cite{ru:VAK:2008.ru}.

%
То, что компонента $R_\iota$ является таким $1$-кодом "--- ее основное свойство.
Поскольку наше определение линейных компонент отличается от других, необходимо доказать этот факт (по существу, следующая лемма совпадает с \cite[Следствие~3.4]{ru:EV:94}).

\begin{lemma}\label{l:O} Для любого $\bar z$ из $F^n$ выполнено тождество
  $\Omega(R_\iota+\bar z) = \Omega(R_\iota+\bar z+\bar e^{(\iota)}).$
\end{lemma}
\proofr Без потери общности положим $\bar z = 0^n$.
Обозначим $\bar e^{(0^m)}:=0^n$. Тогда имеем
$$
\Omega(R_\iota)
=\bigcup_{\kappa\in F^m}(R_\iota + \bar e^{(\kappa)})
=\bigcup_{\kappa\in F^m}(R_\iota + \bar e^{(\iota)} + \bar e^{(\kappa+\iota)})
=\bigcup_{\lambda\in F^m}( (R_\iota + \bar e^{(\iota)}) + \bar e^{(\lambda)})
=\Omega(R_\iota + \bar e^{(\iota)}),
$$
поскольку $\bar e^{(\iota)} + \bar e^{(\kappa)} + \bar e^{(\kappa+\iota)} \in R_\iota$ для всех $\kappa\in F^m$.
\proofend

\begin{lemma}\label{l:RR} Для каждого элемента $\bar c$ множества $<R_\iota,R_\kappa>$ выполнено соотношение
\begin{equation}\label{eq:RR}
 c_{\alpha} + c_{\alpha+\iota}+c_{\alpha+\kappa} + c_{\alpha+\iota+\kappa}=0
\mbox{ для всех } \alpha \in F^m\setminus <\iota,\kappa>
\end{equation}
\end{lemma}
\proofr
По определению, для элементов множеств $R_\iota$ и $R_\kappa$  соотношение (\ref{eq:RR}) выполнено.
Следовательно, это верно и для элементов их линейной оболочки.
\proofend

Следующая лемма является ключевой в наших рассуждениях.

\begin{lemma}\label{l:disj}
  Для любых $\iota,\kappa\in\dot F^m$ на расстоянии не меньше $3$ от $0^m$ и друг от друга, $\iota$"=компонента
  $R_\iota+ \zeroext{\iota} + \bar e^{(\iota)}$
  и
  $\kappa$"=компонента
  $R_\kappa+ \zeroext{\kappa} + \bar e^{(\kappa)}$
  не пересекаются и не содержат $0^n$.
\end{lemma}
\proofr
Исходя из общих алгебраических соображений достаточно показать, что вектор
$\bar w:= \zeroext{\iota} + \bar e^{(\iota)} + \zeroext{\kappa} + \bar e^{(\kappa)}$
 не принадлежит множеству $<R_\iota,R_\kappa>$.
Пусть $j$ "--- ненулевая координата вектора $\iota+\kappa$.
Тогда $\pi^{(j)}$ "--- индекс ненулевой координаты вектора $\bar w$,
индексы других ненулевых координат также принадлежат множеству ${\Pi} \cup \{\iota,\kappa\}$. Но, поскольку попарные расстояния между векторами $0^m$, $\iota$, $\kappa$, и $\iota+\kappa$
не меньше $3$, индексы $\pi^{(j)}+\iota$, $\pi^{(j)}+\kappa$, $\pi^{(j)}+\iota+\kappa$ не принадлежат множеству ${\Pi} \cup \{\iota,\kappa\}$.
Таким образом, имеем
$w_{\pi^{(j)}}+w_{\pi^{(j)}+\iota}+w_{\pi^{(j)}+\kappa}+w_{\pi^{(j)}+\iota+\kappa} = 1+0+0+0=1$
и по лемме~\ref{l:RR} имеем $\bar w \not\in <R_\iota,R_\kappa>$.

Аналогичными рассуждениями доказывается, что ни $R_\iota+ \zeroext{\iota} + \bar e^{(\iota)}$, ни $R_\kappa+ \zeroext{\kappa} + \bar e^{(\kappa)}$ не содержат $0^m$.
\proofend

\begin{theorem}
  Пусть $C\subset F^m$ "--- некоторый $1$-код, содержащий вектор $0^m$; положим $\dot C := C \setminus \{0^m\}$.
Тогда множество
 $$ P(C) := \left( H \bigg\backslash \bigcup_{\iota\in \dot C}(R_\iota+ \zeroext{\iota} + \bar e^{(\iota)})
 \right) \cup \bigcup_{\iota\in \dot C}(R_\iota+ \zeroext{\iota})$$
 является $1$"=совершенным кодом в $F^n$. Более того, верно следующее соотношение:
\begin{equation}\label{eq:CtoP}
  C= \{ \iota \in F^m | (\iota,0^{n-m})\in P(C)\}.
\end{equation}
\end{theorem}
\proofr
Заметим, что согласно (\ref{eq:H}) вектор  $\zeroext{\iota} + \bar e^{(\iota)}$ принадлежит коду $H$ для каждого $\iota$. Таким образом, для каждого $\iota$ выполнено $R_\iota+ \zeroext{\iota} + \bar e^{(\iota)} \subset H$.

По лемме~\ref{l:disj} множества $R_\iota+ \zeroext{\iota} + \bar e^{(\iota)}$, $\iota\in \dot C$,
попарно не пересекаются.
Поскольку они являются подмножествами $1$-совершенного кода, их окрестности также попарно не пересекаются. Учитывая лемму~\ref{l:O}, мы видим, что множество $P(C)$ по определению является $1$-совершенным кодом.

Легко видеть, что\\
(*) \emph{единственное слово в коде $H$, имеющее вид $(\alpha,0^{n-m})$ "--- это слово из всех нулей.}

Далее, $\zeroext{\iota}$ "--- единственное слово такого вида в $R_\iota+ \zeroext{\iota}$, другими словами,\\
(**) \emph{если для некоторого $\kappa\in F^m$
мы имеем $(\kappa,0^{n-m}) \in R_\iota+ \zeroext{\iota}$,
то $\kappa = \iota$.}
В самом деле, предположим, что $\zeroext{\kappa}=(\kappa,0^{n-m})\in R_\iota+ \zeroext{\iota}$.
Тогда $\zeroext{\kappa} + \zeroext{\iota} \in R_\iota \subset H$.
Согласно (*), имеем $\zeroext{\kappa} + \zeroext{\iota}=0^n$,
что доказывает утверждение (**).

Утверждения (*) и (**) доказывают, что соотношение (\ref{eq:CtoP}) вытекает из определения кода $P(C)$.
\proofend

\begin{remark}
Как нетрудно заметить,
ранг (размерности линейной оболочки) построенного кода $P(C)$ вычисляется по формуле ${\mathrm{rank}}(P(C))={\mathrm{rank}}(H)+{\mathrm{rank}}(C)$. Таким образом, начиная с некоторого $m$, мы можем построить код $P(C)$ любого ранга от $2^m-m$ до $2^m-1$.
 Это дает (для больших длин) интересную трактовку результата \cite[Proposition 6.2]{ru:EV:94} о существовании $1$"=совершенных кодов любого ранга от ${\mathrm{rank}}(H)$ до $n$.
\end{remark}
\begin{remark} Как показано в работе \cite{ru:Mal:2001:nonsyst.ru}, достаточно сдвинуть только семь непересекающихся смежных классов $R_\iota+z_\iota$,  $R_\kappa+z_\kappa$,  $R_\lambda+z_\lambda$,
  $R_{\iota+\kappa}+z_{\iota+\kappa}$, $R_{\iota+\lambda}+z_{\iota+\lambda}$,
  $R_{\kappa+\lambda}+z_{\kappa+\lambda}$, $R_{\iota+\kappa+\lambda}+z_{\iota+\kappa+\lambda}$
 в коде $H$ для некоторых линейно независимых $\iota$, $\kappa$, $\lambda$ и некоторых
  $z_\iota,z_\kappa,z_\lambda,z_{\iota+\kappa},z_{\iota+\lambda},z_{\kappa+\lambda},z_{\iota+\kappa+\lambda}\in H$. Таким образом, для любого $m\geq 6$ мы можем построить несистематический $1$"=совершенный код $P(C)$ с $|C|=8$.
\end{remark}

Интересно было бы определить такое минимальное значение величины $n(m)$, что любой $1$-код длины $m$ может быть вложен в $1$"=совершенный код длины $n(m)$. Недавно был решен аналогичный вопрос для систем троек Штейнера \cite{ru:BryHor:Lindner}: для того, чтобы любая частичная система троек Штейнера порядка $u\geq 9$ всегда могла быть вложена в систему троек Штейнера порядка $v$, необходимо и достаточно выполнение условий  $v\geq 2u+1$ и $v \equiv 1,3 \bmod 6$ (последнее условие является необходимым и достаточным для существования систем троек Штейнера порядка $v$). Учитывая связь между совершенными кодами и системами троек Штейнера, упомянутую во введении, мы можем утверждать, что $n(m)\geq 2m+1$ при $m\geq 9$; наша теорема дает оценку сверху $n(m)\leq 2^m-1$.

Авторы благодарят рецензентов за полезные советы.